\documentclass[12pt]{amsart}
\usepackage{amssymb, amstext, amscd, amsmath}
\makeatletter
\def\@cite#1#2{{\m@th\upshape\bfseries%
[{#1\if@tempswa{\m@th\upshape\mdseries, #2}\fi}]}}
\makeatother
\theoremstyle{plain}
\newtheorem{thm}{Theorem}[section]
\newtheorem{cor}[thm]{Corollary}

\newtheorem{lem}[thm]{Lemma}
\theoremstyle{definition}
\newtheorem{rem}[thm]{Remark}
\newtheorem{defn}[thm]{Definition}

\newcommand{\Prf}{\noindent\textbf{Proof.\ }}
\newcommand{\bx}{\hfill$\blacksquare$\medbreak}

\newcommand{\ca}{\mathrm{C}^*}

\newcommand{\cenv}{\mathrm{C}^*_{\text{env}}}

\newcommand{\bbC}{{\mathbb{C}}}

\newcommand{\bbN}{{\mathbb{N}}}

\newcommand{\bbR}{{\mathbb{R}}}
\newcommand{\bbT}{{\mathbb{T}}}
\newcommand{\bbZ}{{\mathbb{Z}}}

  \newcommand{\A}{{\mathcal{A}}}
  \newcommand{\B}{{\mathcal{B}}}

  \newcommand{\F}{{\mathcal{F}}}

  \newcommand{\J}{{\mathcal{J}}}
  \newcommand{\K}{{\mathcal{K}}}
\renewcommand{\L}{{\mathcal{L}}}

\renewcommand{\O}{{\mathcal{O}}}

  \newcommand{\T}{{\mathcal{T}}}

  \newcommand{\X}{{\mathcal{X}}}
  \newcommand{\Y}{{\mathcal{Y}}}

\newcommand{\fT}{{\mathfrak{T}}}

\newcommand{\alg}{\operatorname{alg}}

\begin{document}

%%%%%%%%%%%%%%%%%%%%%%%%%%%%%%%%%%%%%%%%%%%%%%
%%%%%%%%%%%%%%
\title[$\ca$-envelope]{Tensor algebras of $\ca$-correspondences
and their $\ca$-envelopes.}
%\thanks{}
%

\author[E.G. Katsoulis]{Elias~G.~Katsoulis}
%\thanks{Katsoulis' research partially supported by a grant from
%ECU}
\address{Department of Mathematics\\East Carolina University\\
Greenville, NC 27858\\USA}
\email{KatsoulisE@mail.ecu.edu}
\author[D.W. Kribs]{David~W.~Kribs}
\address{Department of Mathematics and Statistics\\University of Guelph\\
Guelph, Ontario\\CANADA N1G 2W1}
\email{dkribs@uoguelph.ca}
\begin{abstract}
We show that the $\ca$-envelope of
the tensor algebra of an arbitrary $\ca$-correspondence $\X$ coincides with the Cuntz-Pimsner
algebra $\O_{\X}$, as defined by Katsura \cite{Ka}. This improves earlier results of Muhly and Solel \cite{MS2} and Fowler,
Muhly and Raeburn \cite{FMR}, who came to the same conclusion under the additional hypothesis that $\X$ is
strict and faithful.
\end{abstract}

\thanks{First author was partially supported by
a summer grant from ECU. Second author was partially supported by an NSERC grant}
\date{}
\maketitle
%%%%%%%%%%%%%%%%%%%%%%%%%%%%%%%%%%%%%%%%%%%%%%
%%%%%%%%%%%%%%

\section{Introduction}\label{S:intro}

Fowler, Muhly and Raeburn have recently
characterized \cite[Theorem 5.3.]{FMR} the $\ca$-envelope of the
tensor algebra $\T^{+}_{\X}$ of a \textit{faithful and strict} $\ca$-correspondence $\X$, as
the associated universal Cuntz-Pimsner algebra. Their proof is based on a
gauge invariant uniqueness theorem and earlier elaborate results of Muhly
and Solel \cite{MS2}. Beyond faithful strict $\ca$-correspondences, little is known: if $\X$ is strict,
but not necessary faithful, then the $\ca$-envelope
of $\T^{+}_{\X}$ is known to be a quotient of the associated
Toeplitz-Cuntz-Pimsner  algebra, without any further
information (Theorem 6.4 in \cite{MS2}). In \cite[Remark 5.4]{FMR}, the authors ask whether
the above mentioned conditions on $\X$ are necessary for the validity of their Theorem 5.3 in
\cite{FMR}.

In this note we answer the question of Fowler, Muhly and Raeburn \cite{FMR} (and Muhly and Solel \cite{MS2}) by showing that the
$\ca$-envelope of
the tensor algebra of an \textit{arbitrary} $\ca$-correspondence $\X$ coincides with the Cuntz-Pimsner
algebra $\O_{\X}$, as defined by Katsura in \cite{Ka}. Our proof does not require any
of the results from \cite{MS2} and is modelled upon the proof
of our recent result \cite{KK2} that identifies
the $\ca$-envelope of the tensor algebra of a directed graph. We also make use of the
result of Muhly and Tomforde \cite{MT} that generalizes the process of adding tails to a graph to
the context of $\ca$-correspondences.

\section{Preliminaries}

Let $\A$ be a $\ca$-algebra and $\X$ be a (right) Hilbert $\A$-module, whose inner product is denoted as
$\langle \, . \mid . \, \rangle$. Let $\L (\X)$ be the adjointable operators
on $\X$ and let $\K(\X)$ be the norm closed subalgebra of $\L(\X )$ generated by the
operators $\theta_{\xi , \eta}$, $\xi , \eta \in \X$, where $\theta_{\xi , \eta}(\zeta)=
\xi \langle \eta | \zeta \rangle$, $\zeta \in \X$.

A Hilbert $\A$-module $\X$ is said to be a \textit{$\ca$-correspondence} over $\A$ provided that there exists
a $*$-homomorphism $\phi_{\X} : \A \rightarrow \L(\X)$. We refer to $\phi_{\X}$ as the left action of $\A$ on
$\X$. A $\ca$-correspondence $\X$ over $\A$ is said to be \textit{faithful}
if and only if the map $\phi_{\X}$ is faithful. A $\ca$-correspondence $\X$ over $\A$ is called \textit{strict} iff
$\overline{[\phi_{\X}(\A)\X ]} \subseteq \X$ is complemented, as a submodule
of the Hilbert $\A$-module $\X$. In particular, if $\overline{[ \phi_{\X}(\A)\X ]} = \X$, i.e., the map
$\phi_{\X}$ is non-degenerate, then $\X$ is said to be \textit{essential}.

From
a given $\ca$-correspondence $\X$ over $\A$, one can form
new $\ca$-correspondences over $\A$, such as the $n$-\textit{fold ampliation} or {direct sum}
$\X^{(n)}$ (\cite[page 5]{L})
and the $n$-\textit{fold interior tensor product} $\X^{\otimes n} \equiv
\X \otimes_{\phi_{\X}}  \X \otimes_{\phi_{\X}} \dots  \otimes_{\phi_{\X}}\X$ (\cite[page 39]{L}, $n \in \bbN$,
($\X^{\otimes 0}\equiv \A$).
These operation are defined within the category of
$\ca$-correspondences over $\A$. (See \cite{L} for more details.)
 
A \textit{representation} $(\pi , t)$ of a $\ca$-correspondence $\X$ over $\A$ on a
$\ca$-algebra $\B$ consists of
a $*$-homomorphism $\pi : \A \rightarrow \B$
and a linear map $t : \X \rightarrow \B$ so that
\begin{itemize}
\item[(i)] $t(\xi )^* t(\eta) = \pi(\langle \xi | \eta \rangle)$, for $\xi , \eta \in \X$,
\item[(ii)] $\pi (a) t(\xi)= t(\phi_{\X}(a) \xi)$, for $a \in \A$, $\xi \in \X$.
\end{itemize}
For a representation $(\pi , t)$ of  a $\ca$-correspondence $\X$ there exists a $*$-homomorphism
$\psi_t : \K(\X ) \rightarrow \B$ so that $\psi_t (\theta_{\xi , \eta}) = t(\xi)t(\eta)^*$,
for $\xi , \eta \in \X$. Following Katsura \cite{Ka}, we say that the representation $(\pi , t)$ is
\textit{covariant} iff $\psi_t ( \phi_{\X} (a)) = \pi(a)$, for all $a \in \J_{\X}$, where
\[
\J_{\X} \equiv \phi_{\X}^{-1}(\K(\X)) \cap (\ker \phi_{\X})^{\perp} .
\]
If $(\pi , t)$ is a representation of $\X$ then the $\ca$-algebra (resp. norm closed algebra) generated by the images of
$\pi$ and $t$ is denoted as $\ca(\pi , t)$ (resp. $\alg((\pi , t)$). There is a
universal representation $(\overline{\pi}_{\A}, \overline{t}_{\X})$
for $\X$ and the
$\ca$-algebra $\ca (\overline{\pi}_{\A}, \overline{t}_{\X})$ is the Toeplitz-Cuntz-Pimsner
algebra $\T_{\X}$. Similarly, the
Cuntz-Pimsner algebra $\O_{\X}$ is the $\ca$-algebra generated by the image of the universal
covariant representation $(\pi_{\A}, t_{\X})$ for $\X$.

A concrete presentation of both $\T_{\X}$ and $\O_{\X}$ can be given in terms of the generalized
Fock space $\F_{\X}$ which we now describe. The \textit{Fock space} $\F_{\X}$ over the
correspondence $\X$ is defined to be the direct sum of the $\X^{\otimes n}$ with the structure of a direct sum
of $\ca$-correspondences over $\A$,
\[
\F_{\X}= \A \bigoplus \X \bigoplus \X^{\otimes 2} \bigoplus \dots .
\]
Given $\xi \in \X$, the (left)
creation operator $t_{\infty}(\xi) \in \L(\F_{\X})$ is defined by the formula
\[
t_{\infty}(\xi)(a , \zeta_{1}, \zeta_{2}, \dots ) = (0, \xi a, \xi \otimes \zeta_1,
\xi \otimes \zeta_2, \dots),
\]
where $\zeta_n \in \X^{\otimes n}$, $n \in \bbN$. Also, for $a \in \A$, we define
$\pi_{\infty}(a)\in \L(\F_{\X})$ to be the diagonal operator with $\phi_{\X}(a)\otimes id_{n-1}$
at its $\X^{\otimes n}$-th entry. It is easy to verify that $( \pi_{\infty}, t_{\infty})$ is a representation of
$\X$ which is called the \textit{Fock representation} of $\X$. Fowler and Raeburn \cite{FR} (resp. Katsura
\cite{Ka}) have shown that
the $\ca$-algebra $\ca ( \pi_{\infty}, t_{\infty})$
(resp $\ca ( \pi_{\infty}, t_{\infty})/ \K(\F_{\X\J_{\X}})$)
is isomorphic to $\T_{\X}$ (resp. $\O_{\X}$).

\begin{defn}
The \textit{tensor algebra} of a $\ca$-correspondence $\X$ over $\A$ is the norm-closed algebra
$\alg(\overline{\pi}_{\A}, \overline{t}_{\X})$ and is denoted as $\T_{\X}^{+}$.
\end{defn}

According to \cite{FR, Ka}, the algebras $\T_{\X}^{+} \equiv \alg(\overline{\pi}_{\A}, \overline{t}_{\X})
$ and $ \alg( \pi_{\infty}, t_{\infty} )$ are completely isometrically isomorphic and we will
therefore identify them. The main result of this paper implies that $\T_{\X}^{+}$ is also
completely isometrically isomorphic to $\alg(\pi_{\A}, t_{\X})$.

\section{Main Result}

We begin with a useful description of the norm in $\X^{(n)}$.

\begin{lem} \label{isomsubst}
Let $\X$, $\Y$ be Hilbert $\A$-modules and let $\phi : \A \rightarrow \L(\Y)$ be an injective
$*$-homomorphism. If $(\xi_i)_{i=1}^{n} \in \X^{(n)}$, then
\begin{equation} \label{norm}
\| (\xi_i)_{i=1}^{n} \| = \sup \{ \| (\xi_i  \otimes_{\phi} u)_{i=1}^{n} \| \mid u \in \Y , \| u \|=1 \}.
\end{equation}

\end{lem}

\Prf Let us denote by $M$ the supremum in (\ref{norm}).
Then, using the fact that $\phi$ is injective and therefore isometric,
\begin{align*}
M^2 &=\sup \{ \| \sum_{i=1}^{n}\,\langle u | \phi(\langle \xi_i |\xi_i \rangle)u \rangle \| \mid u \in \Y , \| u \|=1 \}    \\
 &=\sup \{ \| ( \phi(\langle \xi_i |\xi_i \rangle^{1/2})u)_i \|^{2} \mid u \in \Y , \| u \|=1 \}  \\
  &=\left\|
 \left(
\begin{matrix}
 0 & 0& \ldots &
 \phi(\langle \xi_1 |\xi_1\rangle^{1/2})\cr 0 & 0 &
\ldots & \phi(\langle \xi_2 |\xi_2\rangle^{1/2}) \cr \vdots & \vdots & \ddots & \vdots \cr 0 &
0 & \cdots &\phi(\langle \xi_n |\xi_n \rangle^{1/2})
\end{matrix}
\right) \right\|^2\\
 &=  \| \phi( \sum_{i=1}^{n}\, \langle \xi_i |\xi_i \rangle ) \|  = \| (\xi )_i \|^2
 \end{align*}
and the conclusion follows.
\bx

In the proof of our next lemma we make use of the right creation operators.
If $\Y$ is a $\ca$-correspondence over $\A$
and $\xi \in \Y^{\otimes k}$, then define the right
creation operator $R_{\xi}$ by the formula
\[
R_{\xi}(a , \zeta_{1}, \zeta_{2}, \dots ) = (\underbrace{0, 0, \dots,0}_{k},
(\phi_{\X}(a)\otimes id_{k-1})(\xi),  \zeta_1 \otimes \xi,
  \zeta_2\otimes \xi, \dots),
\]
$\zeta_n \in \Y^{\otimes n}$, $n \in \bbN$. The operator, $R_{\xi}$ may not be adjointable
 but it is nevertheless bounded by $\| \xi \|$ and commutes with $\alg (\pi_{\infty} , t_{\infty})$.

\begin{lem}  \label{fund}
If $\X$ be a faithful $\ca$-correspondence over $\A$, then
\[
\| A \| = \inf \{ \| A  + K \| \mid K \in M_n (\K ( \F_{\X} )) \}
\]
for all $A \in M_n (\T^{+}_{\X })$, $n \in \bbN$.
\end{lem}

\Prf  Let $K \in M_n (\K ( \F_{\X} ))$ be an $n \times n$ matrix with entries in $\K ( \F_{\X} )$
and let $\epsilon >0$. We choose unit vector $\xi \in \ \F_{\X}^{(n)}$ so that
$\| A \xi \| \geq \|A\|-\epsilon $. Since $K \in M_n (\K (  \F_{\X} ))$, there exists $k \in \bbN$
so that $\| K R_{u}^{(n)} \| \leq \epsilon $, for all unit vectors $ u \in \X^{\otimes k}$.
(Here $R_{u}^{(n)}$ denotes the
the $n$-th ampliation of the right creation operator $R_{u}$.)
Note that for any vector $u \in \X^{\otimes k}$ we have
\[
\| R_{u}^{(n)}A\xi \| = \| A\xi \otimes u \|.
\]
Therefore,
using Lemma
\ref{isomsubst}, we choose unit vector
$u \in \X^{\otimes k}$ so that
\begin{align*}
\| R_{u}^{(n)} A\xi \|&\geq \|A\xi \| -\epsilon  \\
                    &\geq \|A \| -2\epsilon .
\end{align*}
We compute,
\begin{align*}
\| A  + K \| &\geq  \| (A  + K)R_{u}^{(n)}\xi \|   \\
             &\geq   \| A R_{u}^{(n)}\xi \| - \epsilon  \\
             &=   \| R_{u}^{(n)}A\xi \| - \epsilon  \\
             &\geq  \|A \| -3\epsilon .
\end{align*}
Since $\epsilon$ and $K$ are arbitrary, the proof is complete.
\bx

\begin{cor}  \label{faithful}
Let $\X$ be a faithful $\ca$-correspondence over $\A$, and let
$(\pi_{\A} , t_{\X})$ be the universal covariant representation of
$\X$. Then, there exists a complete isometry
\[
\tau_{\X}:  \T^{+}_{\X } \longrightarrow \alg (\pi_{\A} , t_{\X})
\]
so that $\tau_{\X}(\pi_{\infty}(a))= \pi_{\A}(a)$, for all $a \in \A$,
and $\tau_{\X}(t_{\infty}(\xi))= t_{\X}(\xi)$, for all $\xi \in \X$.

In particular, the algebra $\alg (\pi_{\A} , t_{\X})$
is completely isometrically isomorphic to the tensor algebra
$\T^{+}_{\X }$.
\end{cor}

\Prf Let $\tau_{\X}$ be the restriction of the natural quotient map
\[
\ca (\pi_{\infty} , t_{\infty}) \longrightarrow \ca (\pi_{\infty} , t_{\infty})  / \K(\F_{\X\J_{\X}})
\]
on the non-selfadjoint subalgebra $\alg (\pi_{\infty} , t_{\infty})$. By Lemma
\ref{fund}, this map is a complete isometry.
\bx

\begin{rem} Note that the above lemma already implies the result
of Fowler, Muhly and Raeburn \cite[Theorem 5.3.]{FMR}
without their requirement of $\X$ being strict.
\end{rem}

We now remove the requirement
of $\X$ being faithful from the statement of the above Lemma. In the special
case of a graph correspondence, this was done in \cite{KK2}
with the help of a well-known process called "adding tails to a graph".
This process has been generalized to arbitrary correspondences
by Muhly and Tomforde \cite{MT}. Indeed, let
$\X$ be an arbitrary $\ca$-correspondence over $\A$ and let $\fT \equiv
c_{0}(\ker \phi_{\X})$ consist of all null sequences in $\ker \phi_{\X}$.
Muhly and Tomforde show that there exists a well defined
left action of $\B \equiv \A \oplus \fT$ on $\Y \equiv \X \oplus \fT$ so that
$\Y$ becomes a \textit{faithful} $\ca$-correspondence over $\B$. One can view $\A$ and the
$\ca$-correspondence $\X$ as a subsets of $\B$ and $\Y$ respectively, via the identifications
\begin{align*}
\A \ni a &\longrightarrow (a, 0) \in \A\oplus 0, \\
\X \ni \xi & \longrightarrow (\xi , 0)\in  \X\oplus 0.
\end{align*}
and by noting that the action of $\phi_{\Y}$ on $\A\oplus 0$ coincides with that of $\phi_{\X}$
on $\A$.
(The restriction of a
representation $(\pi , t) $ of $\Y$ on that subset of $\Y$ will be denoted as
$(\pi_{|\A}, t_{|\X})$ and is indeed a representation of $\X$.) In \cite[Theorem 4.3.(b)]{MT}
it is shown that if $(\pi , t) $ is a covariant representation of $\Y$, then
$(\pi_{|\A}, t_{|\X })$ is a covariant representation of $\X$.

\begin{lem}  \label{ess}
Let $\X$ be a $\ca$-correspondence over $\A$, and let
$(\pi_{\A} , t_{\X})$ be the universal covariant representation of
$\X$. Then, there exists a complete isometry
\[
\tau_{\X}:  \T^{+}_{\X } \longrightarrow \alg (\pi_{\A} , t_{\X})
\]
so that $\tau_{\X}(\pi_{\infty}(a))= \pi_{\A}(a)$, for all $a \in \A$,
and $\tau_{\X}(t_{\infty}(\xi))= t_{\X}(\xi)$, for all $\xi \in \X$.
\end{lem}

\Prf Let $(\pi_{\infty}, t_{\infty})$ be the Fock representation of $\Y$ and note that
\cite[Corollary 4.5]{Ka} shows that
\[
\pi_{\infty}(\B)\, \bigcap \, \psi_{t_{\infty}}(\K(\Y)) = \{ 0 \}.
\]
Therefore,
the restriction $(\pi_{\infty |\, \A}, t_{\infty| \, \X })$ satisfies
the same property and so
\cite[Theorem 6.2]{Ka} implies that the integrated
representation $\pi_{\infty |\, \A} \times t_{\infty| \, \X }$
is a $\ca$-isomorphism from the universal Toeplitz algebra $\T_{\X}$ onto
$\ca (\pi_{\infty |\, \X}, t_{\infty| \, \X })$. We therefore
view $\T^{+}_{\X}$ as a subalgebra of $\T^{+}_{\Y}$.

Corollary \ref{faithful} shows now that there exists a complete isometry
\[
\tau_{\Y}: \T^{+}_{\Y } \longrightarrow \alg (\pi_{\B} , t_{\Y})
\]
so that $\tau_{\Y}(\pi_{\infty}(b))= \pi_{\B}(b)$, for all $b \in \B$,
and $\tau_{\Y}(\pi_{\infty}(\xi))= t_{\Y}(\xi)$, for all $\xi \in \Y$.
As we discussed earlier, \cite[Theorem 4.3.(b)]{MT} shows that the restriction
$(\pi_{\B |\, \A}, t_{\Y | \, \X })$ is covariant for $\X$. Since it is also injective,
the gauge invariant uniqueness theorem \cite[Theorem 6.4]{Ka} shows that the restriction
$\tau_{\X} \equiv \tau_{\Y}|_{\T_{\X}}$ has range isomorphic to $\alg (\pi_{\A} , t_{\X})$
and satisfies the desired properties.
\bx

Let $\B$ be a $\ca$-algebra and let $\B^{+}$ be a (nonselfadjoint)
subalgebra of $\B$ which generates $\B$ as a $\ca$-algebra and
contains a two-sided contractive approximate unit for $\B$, i.e.,
$\B^{+}$ is an essential subalgebra for  $\B$. A two-sided ideal $\J$
of $\B^{+}$ is said to be a \textit{boundary ideal} for $\B^{+}$ if and
only if the quotient map $\pi: \B \rightarrow \B / \J$ is a
complete isometry when restricted to $\B^{+}$. It is a result of
Hamana \cite{Hamana}, following the seminal work of Arveson
\cite{Arvenv}, that there exists a boundary ideal $\J_S (\B^{+})$, the
\textit{Shilov boundary ideal}, that contains all other boundary
ideals. In that case, the quotient $\B / \J_S (\B^{+})$ is called the
\textit{$\ca$-envelope} of $\B^{+}$ and it is denoted as $\cenv(\B^{+})$. The $\ca$-envelope is unique in
the following sense: Assume that $\phi': \B^{+} \rightarrow \B'$ is
a completely isometric isomorphism of $\B^{+}$ onto an essential
subalgebra of a $\ca$-algebra $\B'$ and suppose that the Shilov
boundary for $\phi'(\B^{+}) \subseteq \B'$ is zero. Then $\B$ and
$\B'$ are $*$-isomorphic, via an isomorphism $\phi$ so that
$\phi( \pi(x)) = \phi'(x)$, for all $x \in \B$.

In the case where an operator algebra $\B^{+}$ has no contractive approximate identity,
the $\cenv(\B^{+})$ is defined by utilizing the unitization \cite{Mey} $(\B^{+})_{1}$ of $\B^{+}$:
the $\ca$-envelope of $\B^+$ is the $\ca$-subalgebra of $\cenv((\B^{+})_{1})$ generated by $\B^{+}$.
(See \cite{BL, BS} for a comprehensive discussion regarding the implications of \cite{Mey}
on the theory of $\ca$-envelopes.)

\begin{lem} \label{unitiz}
Let $\B$ be a non-unital $\ca$-algebra and let $\J \subseteq \B_1$ be a closed two-sided
ideal in its unitization. If $\J \cap \B = \{ 0\}$ then $\J = \{ 0\}$.
\end{lem}

\Prf Assume that $\J \neq \{0 \}$. Since $\B_1 \subseteq \B$ has codimension $1$, $\J$ is of the form $\J= [\{ B +\lambda I \}] $, for some $B \in \B$ and non-zero $\lambda \in \bbC$.
Then, easy manipulations show that there is no loss of generality assuming
that $\lambda \in \bbR$ (because $\J \J^* \neq 0$), $A$ is selfadjoint
(because $\J\cap \J^* \neq 0$) and
\begin{equation} \label{trivial}
(A + \lambda I)^2 = A + \lambda  ,
\end{equation}
after perhaps scaling (since $\J^2 \neq 0$). It is easy to see
now that (\ref{trivial}) implies that $A = -P$, for some projection $P\in \B$. But then, $(I-P)\B = 0$ and so
$P$ is a unit for
$\B$, a contradiction.
\bx

We have arrived to the main result of the paper.

\begin{thm} \label{main}
 If $\X$ is a $\ca$-correspondence over $\A$, then the $\ca$-envelope of
$\T^{+}_{\X}$ coincides with the universal Cuntz-Pimsner algebra $\O_{\X}$.
\end{thm}

\Prf According to Lemma \ref{ess}, it suffices to show that
the $\ca$-envelope of $\alg (\pi_{\A} , t_{\X})$ equals $\O_{\X}$.

Assume first that $\alg (\pi_{\A} , t_{\X})$ is unital. In
light of the above discussion,
we need to verify that the Shilov
boundary ideal $\J_S(\alg (\pi_{\A} , t_{\X}))$ is zero.
However, the maximality of $\J_S(\alg (\pi_{\A} , t_{\X}))$ and the invariance
of $\alg (\pi_{\A} , t_{\X})$ under the gauge action of $\bbT$ on $\O_{\X}$ imply
that $\J_S(\alg (\pi_{\A} , t_{\X})\,)$ is a gauge-invariant ideal. By the gauge invariant uniqueness theorem
\cite[Theorem 6.4]{Ka},
any non-zero gauge-invariant ideal
has non-zero intersection with $\pi_{\A}(\A)$. Hence
$\J_S(\alg (\pi_{\A} , t_{\X})\,) = \{ 0 \}$, or otherwise the quotient
map would
not be faithful on $\alg (\pi_{\A} , t_{\X})$.

Assume now that $\alg (\pi_{\A} , t_{\X})$ is not unital. We distinguish two cases.

If $\O_{\X}$ has a unit $I \in \O_{\X}$ then let
\[
\alg (\pi_{\A} , t_{\X})_1 \equiv \alg (\pi_{\A} , t_{\X}) + \bbC I \subseteq \O_{\X}.
\]
Clearly, $\alg (\pi_{\A} , t_{\X})_1$ is gauge invariant and so a repetition of
the arguments in the second paragraph of the proof shows that
\[\cenv(
\alg (\pi_{\A} , t_{\X})_1) =\O_{\X}.
\]
The $\ca$-subalgebra of $\O_{\X}$
generated by  $\alg (\pi_{\A} , t_{\X})$ equals $\O_{\X}$, which by convention will be its
$\ca$-envelope.

Finally, if  $\O_{\X}$ does not have a unit then unitize $\O_{\X}$ by joining a
unit $I$ and let
\[
\alg (\pi_{\A} , t_{\X})_1 \equiv \alg (\pi_{\A} , t_{\X}) + \bbC I
\subseteq \O_{\X} +\bbC I .
\]
Since the Shilov ideal $\J_S (\alg (\pi_{\A} , t_{\X})_1)$ is gauge invariant,
\[
\J_S (\alg (\pi_{\A} , t_{\X})_1) \cap \O_{\X} \subseteq \O_{\X}
\]
is gauge invariant. Therefore,
\[
\J_S (\alg (\pi_{\A} , t_{\X})_1) \cap \O_{\X} =\{0\} ,
\]
or else it meets $\pi_{\A}(\A)$.
By Lemma \ref{unitiz}, $\J_S (\alg (\pi_{\A} , t_{\X})_1) = \{0\}$ and so
$\cenv(\alg (\pi_{\A} , t_{\X})_1) =\O_{\X} +\bbC I$. The $\ca$-subalgebra of $\O_{\X} +\bbC I$
generated by  $\alg (\pi_{\A} , t_{\X})$ is $\O_{\X}$, and the conclusion follows.
\bx

\begin{rem} In \cite[page 596]{FMR}, it is claimed that if a $\X$ is a $\ca$-correspondence
over $\A$, with universal Toeplitz representation $( \overline{\pi}_{\A}, \overline{t}_{\X})$,
then $\overline{\pi}_{\A}$ maps an approximate unit of $\A$ to an approximate unit for both
$\T_{\X}$ and $\T_{\X}^{+}$. It is not hard to see that this claim is valid if and only if $\phi_{\X}$
is non-degenerate. Therefore, there is a gap in the proof of \cite[Theorem 5.3]{FMR} in the case where
$\X$ is strict but not essential. Nevertheless, our Theorem \ref{main}
incorporates all possible cases and hence completes the proof of \cite[Theorem 5.3]{FMR}.
\end{rem}

We now obtain one of the main results of \cite{KK2} as a corollary.

\begin {cor}{\cite[Theorem 2.5]{KK2}}  \label{genvelope}
If $G$ is a countable directed graph then the $\ca$-envelope of $\T_{+}(G)$ coincides with the
universal Cuntz-Krieger algebra associated with $G$.
\end{cor}

Note that in \cite{KK2}, the proof of the above corollary is essentially self-contained and
avoids the heavy machinery used in this paper. The reader would actually benefit from reading that proof
and then making comparisons with the proof of Theorem \ref{main} here.

%%%%%%%%%%%%%%%%%%%%%%%%%%%%%%%%%%%%%%%%%%%%%%%%%%%%%%%%%%%%%%%%%%%5

\vspace{0.1in}

{\noindent}{\it Acknowledgement.} We are grateful to David Blecher for directing us to the work
of Meyer \cite{Mey} and its impact on the theory of $\ca$-envelopes \cite{BL}.

%%%%%%%%%%%%%%%%%%%%%%%%%%%%%%%%%%%%%%%%%%%%%%%%%%%%%%%%%%%

%%%%%%%%%%%%%%%%%%%%%%%%%%%%%%%%%%%%%%%%%%%%%%%%%%%%%%%%%%

\end{document}